\documentclass{article}

\usepackage{arxiv}

\usepackage[utf8]{inputenc} 
\usepackage[T1]{fontenc}    
\usepackage{hyperref}       
\usepackage{url}            
\usepackage{booktabs}       
\usepackage{amsfonts}       
\usepackage{nicefrac}       
\usepackage{microtype}      
\usepackage{lipsum}		
\usepackage{graphicx}
\usepackage{doi}

\usepackage[english]{babel}
\usepackage{amsmath}
\usepackage{epstopdf}
\usepackage{flushend}
\usepackage{float}
\usepackage{graphics}

\usepackage{amssymb}
\usepackage[figuresright]{rotating}
\usepackage{subfigure}

\usepackage{adjustbox}
\usepackage{array}
\makeatletter
\newcommand{\thickhline}{%
    \noalign {\ifnum 0=`}\fi \hrule height 1pt
    \futurelet \reserved@a \@xhline
}
\newcolumntype{"}{@{\hskip\tabcolsep\vrule width 1pt\hskip\tabcolsep}}
\makeatother

\usepackage{multirow}

\usepackage{xcolor}

\DeclareMathOperator*{\argmin}{arg\,min}
\usepackage{algorithm}
\usepackage{algpseudocode}
    \makeatletter
    \def\BState{\State\hskip-\ALG@thistlm}
    \makeatother

\title{RISING a new framework for few-view tomographic image reconstruction with deep learning}

\author{ Davide Evangelista \\ 
	Department of Mathematics\\
	University of Bologna, Italy\\
	\texttt{davide.evangelista5@unibo.it} \\
	\And
	 Elena Morotti \\ 
	Department of Political and Social Sciences\\
	University of Bologna, Italy \\
	\texttt{elena.morotti4@unibo.it}
	\And
	 Elena Loli Piccolomini \\ 
	Department of Computer Science and Engineering\\
	University of Bologna, Italy \\
	\texttt{elena.loli@unibo.it} \\
}

\renewcommand{\shorttitle}{RISING a new framework for few-view tomographic image reconstruction with DL}

\hypersetup{
pdftitle={A template for the arxiv style},
pdfsubject={q-bio.NC, q-bio.QM},
pdfauthor={David S.~Hippocampus, Elias D.~Striatum},
pdfkeywords={First keyword, Second keyword, More},
}

\begin{document}
\maketitle

\begin{abstract}
This paper proposes a new two-step procedure for sparse-view tomographic image reconstruction. It is called RISING, since it combines an early-stopped Rapid Iterative Solver with a subsequent Iteration Network-based Gaining step.
So far, regularized iterative methods have widely been  used for X-ray computed tomography image reconstruction from low-sampled data, since they converge to a sparse solution  in a suitable domain, as upheld by the Compressed Sensing theory. Unfortunately, their use is practically limited by their high computational cost which imposes to perform only a few iterations in the available time for clinical exams. 
Data-driven methods, using neural networks to post-process a coarse and noisy image obtained from geometrical algorithms, have been recently studied and appreciated for both their computational speed and accurate reconstructions. However, there is no evidence, neither theoretically nor numerically, that neural networks based algorithms  solve the mathematical inverse problem modelling the tomographic reconstruction process. 
In our two-step approach,  the first phase  executes very few iterations of a regularized model-based algorithm whereas the second step completes the missing iterations by means of a neural network.
The resulting hybrid deep-variational framework preserves the convergence properties of the iterative method and,  at the same time, it exploits the computational speed and flexibility of a data-driven approach.
Experiments performed on a simulated and a real data set confirm the numerical and  visual accuracy of the reconstructed RISING images in short computational times.
\end{abstract}


\section{Introduction \label{intro}}

Combining healthy protocols with high quality  images  is one of the most important component of medical imaging and a crucial target for researchers involved in minimal invasive Computed Tomography (CT).    
Radiologists, manufacturers and medical physicists have implemented many examination protocols as well as software and hardware modifications to reduce the  harmful ionizing radiations and pave the way to X-ray examinations for screening tests, pediatric cases, or pre-surgical examinations. 
There are two main techniques allowing for a significant reduction of the total  radiation exposure per patient. 
The first one consists in reducing the X-ray tube current  at each scan, without changing the full geometry traditionally used in CT  ({\it low-dose CT}). The measurement data are  very noisy  due to the excessive quantum noise.
The second practical way to lower the radiation per person  consists in reducing the number of X-ray projections ({\it few-view CT}), which leads to incomplete tomographic data, but very fast examinations.
In this paper, we focus on few-view CT images, whose reconstruction is tricky: 
the lack of projections makes the discrete inverse problem mathematically modelling the reconstruction phase have infinite possible solutions \cite{mueller2012linear}.
Conventional filtered backprojection (FBP) algorithms, widely exploited in classical CT, do not  provide stable reconstructions as the Tuy-Smith condition is not satisfied \cite{tuy1983inversion}, and the recovered images suffer from severe striking artifacts. \\ 
A widely used alternative approach is represented by model-based iterative  methods. They model the image reconstruction as a mathematical linear inverse problem which is solved, in the discrete setting, by minimizing a constrained or unconstrained function combining a data-fit term and a regularizer.
The embedding of a sparsifying regularizer  mitigates the lack of many angled views, according to the Compressed Sensing (CS) theory \cite{donoho2006compressed}. However, the resulting iterative schemes solving the  minimization problem are still computationally expensive and typically need several iterations to achieve high quality results.
 An exhaustive review of model-based reconstrction methods can be found in \cite{GraffSidky2016}.

Recently, Deep Learning (DL) based methods have emerged over fully conventional or variational approaches for few-view tomographic reconstruction \cite{wang2018image}. 
The two widely used strategies consist in  unrolled methods, which mimic the implementation of an iterative reconstruction algorithm through network layers, and the so called Learnt Post Processing (LPP) approach.
In particular, the LPP is a two-step scheme where first a low quality image with artifacts and noise is reconstructed with a fast method (typically a FBP)  and then a neural network  suppresses the artifacts. Usually, the network  learns from a set of {\em ground truth} images reconstructed from full dose acquisitions. 
The pioneering works by Han demonstrated in 2016 and 2018 the superiority of LPP strategies over some model-based iterative algorithms for sparse-view CT images \cite{han2016deep,  han2018framing}. 
However, in their inspiring work \cite{sidky2020cnns}, Sidky et al. claim that the popular LPP schemes lack of mathematical characterization.

\begin{figure*}[ht]
   \centering
	\includegraphics[width=0.7\textwidth]{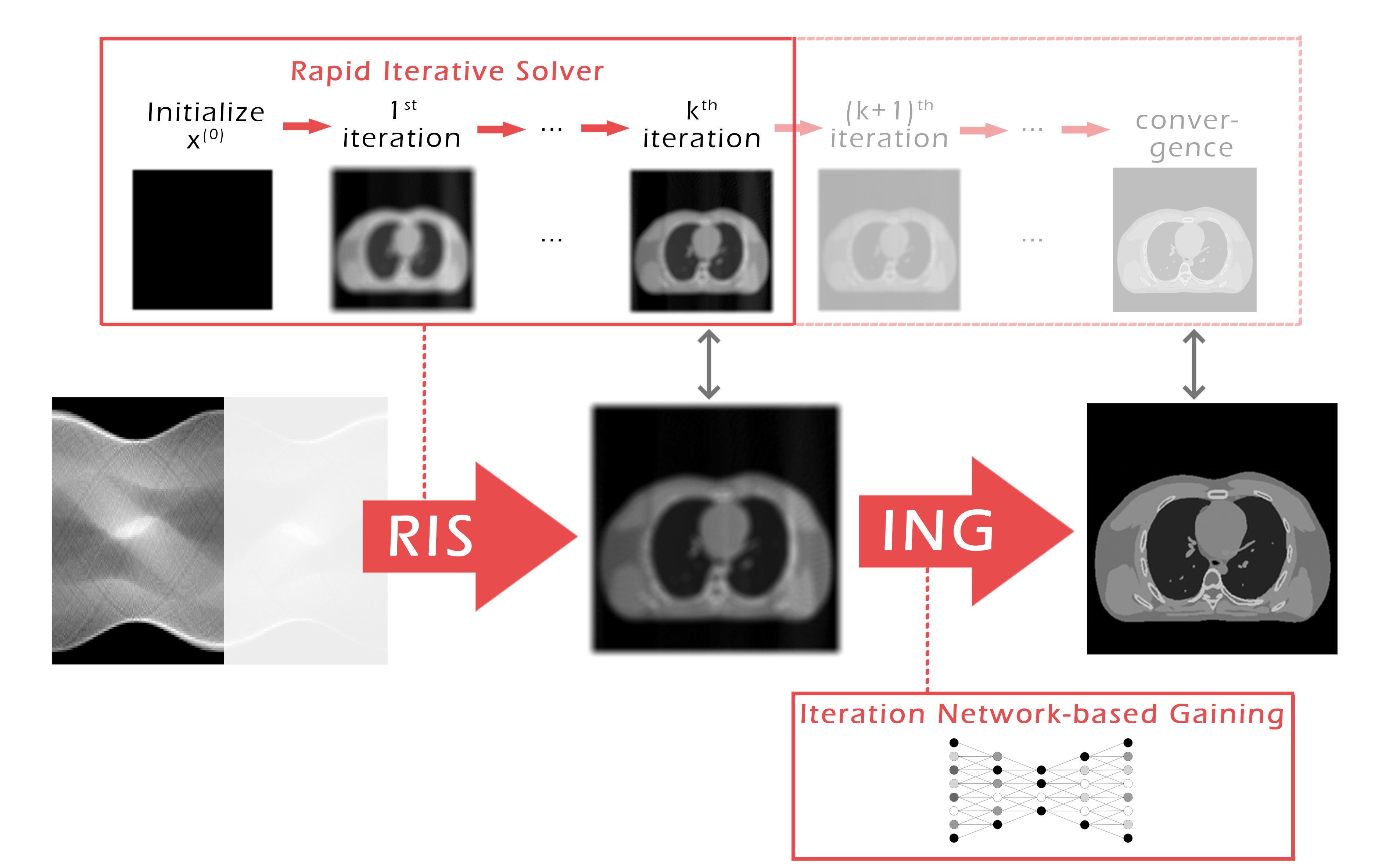}
    \caption{Graphical draft of the proposed two-step RISING workflow for tomographic reconstruction from sparse-view data.}
    \label{fig:RISING}
\end{figure*}

{\bf Motivations.}
The main disadvantage of model-based iterative reconstruction algorithms is their high computational cost. 
On real systems, in order to fulfill the clinical requirements of time per exams, only very few iterations of the algorithms can be performed, leading to a solution which is not the best one. However, by considering fast iterative methods, the images reconstructed in few  iterations already contains many details of interest of the scanned object \cite{Scirep20, piccolomini2021model}.\\
Focusing on data-driven methods based on deep learning, the work by Sidky et al. \cite{sidky2020cnns} claims and demonstrates that LPP schemes do not compute the solution of the CT inverse problem and they can introduce structures not belonging to the scanned objects in the reconstructions . On the contrary, numerical evidence shows that model-based methods compute a good solution of the inverse problem, in the sense that will be defined in the next section. 
In addition, it is worth notice that the efficiency of the deep learning  based methods strictly depends on the training phase and, above all for CT applications, on the training samples. 
In fact, neural networks  need to be trained on task-specific data sets to properly learn both the degradation effects to remove and the anatomical details to preserve characterizing each typology of medical imaging, and a further relevant challenge for medical applications is the lack of precise training data \cite{wang2018image}.
To the best of our knowledge, in the state-of-art literature the neural networks are always trained on ground truth images, achieved with full-dose CT protocols, hence we argue that such projection measurements may be not achievable in some CT applications (for instance, where physical constraints prevents the full-range acquisition, as in limited-angle protocols).


{\bf Contribution.}
Strong of this awareness, this paper proposes a new hybrid iterative and data-driven scheme for the resolution of the sparse-sampling CT inverse problem.

This work has a dual purpose. On one side, we aim at computing a suitable solution of  the inverse  tomographic reconstruction problem as accurately as possible in a short computational time, consistent with the clinical requirements.
On the other side, we intend to remedy to the lack of task-specific CT data sets by using, for networks training, images created by the same system and under the same geometry used for the reconstruction. This is attainable for every CT system.

To achieve our purposes, we propose to apply a model-based iterative method and speed up the approaching to convergence with a neural network trained on accurate task-specific images achieved by few-view protocols.

\textbf{Proposal.}
We refer to our proposal as the {\it RISING} (Rapid Iterative Solver with  Iteration Network-based Gaining) framework.
RISING is conceived as the reconstructing procedure, described by the two following steps executed in sequence:
\begin{itemize}
    \item  
     Staring from the subsampled projection data, a rapid iterative algorithm produces a preliminary coarse reconstruction, by solving  the model-based problem handling few-view CT with a few iterations. The  execution of only few iterations  fits realistic time constraints. 
    
    \item 
    The previously computed rough reconstruction is processed by a deep neural network which aims at retrieving the unperformed iterations towards the convergence. Its output is the RISING solution image.
    As the execution of a neural network is very fast, 
    the use of DL in this step greatly speed-ups the whole reconstruction.
    Differently from other DL-based schemes, in RISING the network is trained on target images obtained offline  as solutions of iterative solver at convergence, from projection data acquired under the CT geometry in exam. This ensures a fully consistent training to all types of medical images and sparse-view protocols.
\end{itemize}   
A graphical draft of RISING is depicted in Figure \ref{fig:RISING}. \\
To validate our proposal, we apply the RISING approach on real and synthetic CT data sets simulating  various sparse-view settings.
In all the cases, we achieve high quality reconstructions and the quality measures on the synthetic data set confirm that the results are reliable solutions of the CT inverse problem. We developed our RISING workflow in Python and the code is available at \url{https://github.com/loibo/RISING}.

\textbf{Organization of the paper.} 
The paper is organized as follows. In Section \ref{sec:recons} we state the CT inverse problem and present the state-of-art DL-based methods addressing it, then Section \ref{sec:rising}  describes the proposed RISING framework.
In Section \ref{sec:ImplemNotes} we introduce the experimental settings considered to achieve the results reported and discussed in Section \ref{sec:results}. 
At last, final conclusions are drawn in Section \ref{sec:concl}.



\section{Few-view X-ray CT image reconstruction methods \label{sec:recons}}

It is well known that the tomographic acquisition  can be modelled as a linear process where the projections are obtained by applying the Radon transform to the object \cite{kak2001principles}.



For simplicity, we consider in this paper the model discretization in two dimensions. The extension to three dimensions is straightforward.\\
In  discrete setting,  the CT process of X-ray absorption is expressed as:
\begin{equation}\label{eq:LinearSyst}
    Ax = b
\end{equation}
where the unknown vector $x \in \mathbb{R}^{n^2}$ denotes the image to reconstruct, the right-hand side term $b$ is an $m$-dimensional vector containing the projection measurements and the $m \times n^2$ system matrix $A$ is the discretization of  the X-ray physical process projecting an $n\times n$ image onto the detector.
We say that an inverse problem is {\it solvable} in a subset $\Omega$ if it admits  a unique solution in $\Omega$.

In case of few-view CT protocols, equation \eqref{eq:LinearSyst} is not solvable in $\mathbb{R}^{n^2}$, since  the under-determined linear system admits infinite solutions (recalling that $m<n^2$). 
According to the Compressed Sensing (CS) theory \cite{donoho2006compressed}, if the desired solution $x^*$ of \eqref{eq:LinearSyst} is sparse in some transform  $Tx$, then \eqref{eq:LinearSyst} is solvable in the subset $\Omega=\{x \in \mathbb{R}^{n^2} | \ ||Tx||_0 \le d \}$, where $\| \cdot \|_0$ is the $\ell_0$ semi-norm counting the non-zero elements of the vector argument.
 Hence the CT inverse problem can be reformulated as the following minimization:

\begin{equation}
x^* = \argmin_x \|Tx\|_0 \text{ s.t. } Ax=b.
\label{eq:L0_CS}
\end{equation}

Since equation  \eqref{eq:L0_CS} is difficult  to solve numerically, it is usually relaxed as:
\begin{equation}
\argmin_x \|Tx\|_1 \text{ s.t. } Ax=b,
\label{eq:4}
\end{equation}
where $\| \cdot \|_1$ is the $\ell_1$ norm \cite{candes2006robust}. 
An unconstrained formulation of \eqref{eq:4} can be stated as:
\begin{equation} 
    \argmin_{x} \| Ax-b\|_2^2 + \lambda \|Tx\|_1
    \label{eq:6}
\end{equation}
where $\lambda$ is a suitable positive parameter (also called {\em regularization parameter}).
The resulting convex optimization problem can be solved by convergent iterative methods such as the lagged diffusivity fixed point algorithm \cite{piccolomini2016fast} or FISTA \cite{xu2016accelerated}. 
Usually,  the  constraint $x \geq 0$ is added in \eqref{eq:4} or \eqref{eq:6}, to preserve the physical the non-negativity properties of the  attenuation coefficients \cite{sidky2014constrained,  purisha2017controlled, piccolomini2021model}.

\subsection*{Deep methods for CT reconstruction \label{sub:deep}}
As anticipated, algorithms exploiting deep Convolutional Neural Networks (CNNs) have been proposed for few-view CT image reconstruction for few years.
Many publications have focused on learning iterative schemes through {\it unrolling} (or unfolding) strategies, where a model-based iterative solver is unfolded into the sequence of its iterations \cite{monga2021algorithm}.
Since each iteration of the unrolled methods matches the  computations performed in the correspondent iteration of the model-based algorithm, these methods aim at recovering the solution of the inverse problem through a minimization. However, they are trained on full dose images as targets, which differ from the inverse problem solution. Hence, there is no evidence that the reconstructed image well approximates the solution of the inverse tomographic problem.
Moreover, the  neural networks employed for unfolded schemes suffer from instabilities, such as vanishing gradient \cite{bengio_book}. To reduce these unwanted effects only a limited number of unfolded iterations are performed, stopping the method far before convergence.\\
The proposals differ for the considered iterative scheme and for the block-per-iteration learnt by the neural network. We give a short overview of the most recent works.
In \cite{adler2017solving}, Adler and Oktem  have developed a partially learned gradient descent algorithm, whereas they have worked on the Chambolle-Pock scheme in \cite{adler2018learned}. 
In \cite{gupta2018cnn} a CNN is trained to act like a projector in a gradient descent algorithms, while in \cite{fista_net}, both the proximal operator and gradient operator (which correspond to $A^T$ when it is applied to a problem corrupted with additive Gaussian noise) of an unrolled FISTA scheme are learned, whereas 
in \cite{zhang2020metainv} the neural network learns the initial iterate of the inner Conjugate Gradients solver in a splitting scheme for optimization.

A completely different data-driven proposal for CT image reconstruction is represented by the {\it Learnt Post-Processing} (LPP) approach. 
It consists in a two-step framework where, firstly, a fast reconstruction algorithm is applied  off-line to compute the medical image from the projection data, then a learnt neural network is applied as post-processing operator to restore the low-quality image. 
It is worth notice  that  most of the networks used for LPP are applied to corrupted images provided by the FBP solver (see for example  \cite{han2018framing,pelt2018improving,zhang2019dualres,schnurr2019simulation,urase2020simulation,morotti2021green} and the references therein).  
Only in \cite{bubba2019learning} and in \cite{jiang2019augmentation} the preliminary solutions are computed by wavelet- or TV-based iterative solvers, respectively.

{\bf Model-based vs LPP methods.}
In addition, the seminal paper \cite{sidky2020cnns} has found out an important observation addressing  LPP schemes in comparison to CS model-based approaches.
Indeed, CS theory has not provided analytic results for X-ray CT measures, hence the authors have considered a set of digital breast phantoms and built CT simulations. 
When considering as ground truth  test images characterized by sparse gradient domains, 
the tests confirmed that the solution of problem \eqref{eq:4} computed by a regularized CS based iterative method is very close to the true solution.
On the contrary, the LPP framework produces worse performances, showing that the  network is not able to solve the  inverse problem. 

\section{The proposed RISING framework \label{sec:rising}} 


In this section, we describe the proposed new framework. 
It is denoted as RISING since it combines a Rapid Iterative Solver (RIS) phase with a following Iteration Network-based Gaining (ING) phase. 
The RISING workflow is also graphically depicted in Figure \ref{fig:RISING}.

Briefly, we look for a sparse solution in a gradient domain and we apply an iterative model-based algorithm to a minimization model such as \eqref{eq:4}.
Instead of getting the algorithm solution at convergence,  we perform only a predefined (and relatively small) number $K$ of iterations and stop the solving algorithm before its convergence, to meet the constraints imposed by clinical setting where very short computational time are admitted to the reconstruction process. 
The resulting early output represents the input of the neural network which enhances the image in a new perspective, trying to retrieve the unperformed iterations instead of executing a LPP.




\subsection{Rapid Iterative Solver  \label{subsec:RIS}} 

The first step of the RISING framework implements the Rapid Iterative Solver (RIS) to the following  minimization problem:
\begin{equation} 
    \argmin_{x \geq 0} \| Ax-b\|_2^2 + \lambda TV_{\beta}(x)
    \label{eq:7}
\end{equation}
where, in order to have a differentiable objective function, we consider the smoothed version of the Total Variation (TV) \cite{RUDIN1992} operator defined as:
\begin{equation} \label{eq:TV}
    TV_{\beta}(x) = \sum_{j=1}^{n^2} \sqrt{ \| \nabla x_j \|_2^2 + \beta^2}
\end{equation}
with $\beta$  fixed  small positive parameter.


Among the wide class of iterative solvers for the problem  \eqref{eq:7}, we select the Scaled Gradient Projection  (SGP) algorithm, proposed in 2008 for image deblurring  \cite{bonettini2008scaled}. The SGP has been  successfully applied with some acceleration techniques for few-view CT reconstructions   \cite{coap2018, Scirep20, piccolomini2021model}, where the authors show t the objects of interest are distinguishable in the solution obtained after very few iterations of the method. 
Regarding the convergence, it is proved in 
 \cite{BPR16} that the theoretical convergence rate  of the SGP at the unique minimum of \eqref{eq:6}  is $\cal{O}$(1/k); they also  numerically show that the practical performance of SGP method  is very well comparable with the convergence rate of the optimal algorithms.
For a detailed presentation of SGP method in tomographic imaging see \cite{coap2018}.

In RISING, we stop the  execution of the  Iterative Solver at a prefixed iteration $K$, far before the numerical convergence, and  the  early solution $x^{(K)}$ is denoted as $x_{RIS}$ in the following.
From our previous works \cite{coap2018, Scirep20} we know that,
due to the empirical rapidity of the SGP, in the achieved coarse image reconstruction the anatomical structures are present but typically very blurred and weakly evident.

\subsection{Iteration Network-based Gaining  \label{subsec:ING}}

The second step of the RISING framework implements the Iteration Network-based Gaining (ING) task. 
Here, a Convolutional Neural Network (CNN) learns the transformation mapping the early solution $x_{RIS}$, achieved in $K$ iterations, to the corresponding convergence image computed in $K^*$ iterations. 

We denote with $\mathcal{T}_j$ the function describing the action of the $j$-th iteration of the solver such as:
$$
 x^{(j+1)} = \mathcal{T}_j(x^{(j)};b) \quad \forall j \geq 0.
$$
Thus, the whole iterative process can be expressed as the concatenation of the following $g_K$ and $f_K$ functions:
\begin{equation}
    g_K := \mathcal{T}_{K-1} \circ \dots \circ \mathcal{T}_{1} \circ \mathcal{T}_{0}
    \label{eq:gComposition}
\end{equation}
and 
\begin{equation}
    f_K := \mathcal{T}_{K^*} \circ \dots \circ \mathcal{T}_{K+1} \circ \mathcal{T}_{K}.
    \label{eq:fComposition}
\end{equation}
In fact, if we label the image achieved by the Iterative Solver at convergence as $x_{IS}$, it holds:
\begin{equation}
   f_K(g_K(x^{(0)};b))= x_{IS}. 
\end{equation}
Since in the RIS step we compute $x_{RIS} = g_K(x^{(0)};b)$, in the ING phase we need to recover $f_K$. 
To do so, we train a CNN to play as $f_K$ when applied on $x_{RIS}$.
Ideally, if the network would perfectly learn $f_K$, its output $x_{ING}$ should be equal to $x_{IS}$. 
It is well known that this is not doable in practice; nevertheless, we show numerical evidence of our claims, through very accurate RISING reconstructions on simulations, in section \ref{sec:results}.

We remark that, differently from the RISING approach,  in case of LPP the CNN tries to learn an image-to-image relation which is not necessarily a function.

\section{Experimental design and implementation notes \label{sec:ImplemNotes}} 

To numerically verify the feasibility of the RISING workflow, we develop numerical simulations.
We consider real medical images, to properly understand the potential of RISING for clinical applications and settings, as well as a synthetic data set of images with geometric elements, where measures of merits can be computed.

This section describes the experimental design, whereas the results are reported and discussed in section \ref{sec:results}.

\subsection{Data set and test problems  from real medical images \label{sub:mayo}}

\begin{figure}[ht]
    \centering
    \includegraphics[height=3.5cm]{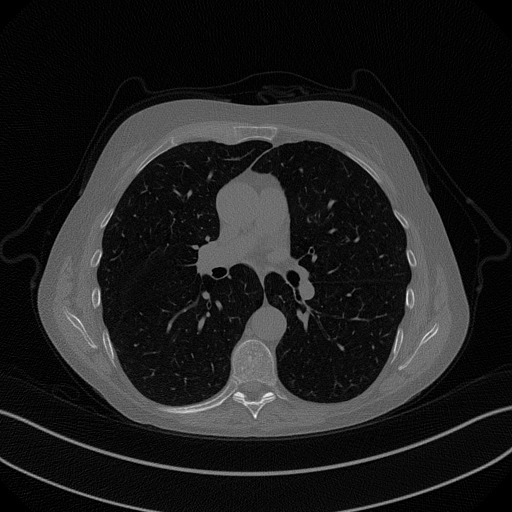}
    \includegraphics[trim= 190 122 130 160, clip, height=3.5cm]{imm/GT.png}
    \includegraphics[trim= 282 102 27 160, clip, height=3.5cm]{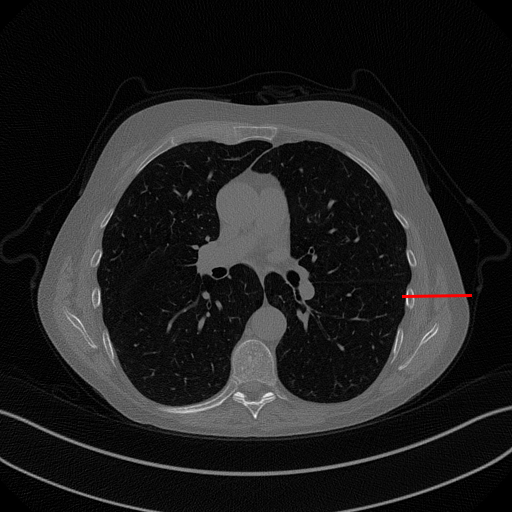}
    \caption{A ground-truth image from the Low Dose Mayo data set, with two zoomed crops on regions with different anatomical structures.}
    \label{fig:gtMAyo}
\end{figure}

As real patient images, we have downloaded the widely used AAPM Low Dose CT Grand Challenge data set by the Mayo Clinic \cite{mccollough2016tu}.
The considered images are $512 \times 512$ pixel reconstructions of human abdomen, computed from full-dose acquisitions. In Figure \ref{fig:gtMAyo} we depict one image with two zooms-in highlighting areas with different anatomical structures, such as pulmonary details, sections of ribs and low-contrast inter-costal muscles.
In all the experiments reported in \ref{ssec:ResultsMayo}, we have used the images from the data set as ground truth $x_{GT}$ references. 
Coherently, we simulate the tomographic projections of the ground truth images, according to a 2D fan-beam geometry, and we add to the sinograms white Gaussian noise with $10^{-2}$ noise level.
To address sparse-view CT reconstructions, we considered two different protocols: the first one is a full angular acquisition with 1-degree spaced projections (we call it $P_{360, 360}$ in the following); in the second one the scanning trajectory  covers $180$ degrees and computed only $60$ projections (it is labelled as $P_{180, 60}$).\\
The RIS step executes only $K = 15$ iterations, hence we denote as $x_{RIS} = x^{(15)}$ the SGP output which is passed as input to the ING phase.


\subsection{Data sets and test problems from synthetic images  \label{ssec:COULE}}

\begin{figure}[ht]
    \centering
    \includegraphics[height=3.5cm]{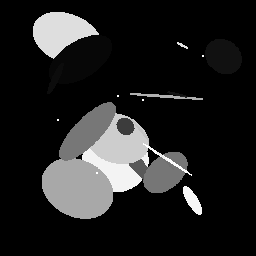}
    \includegraphics[height=3.5cm]{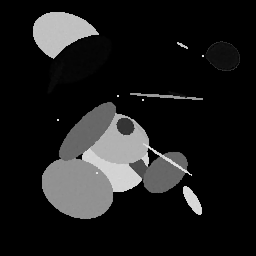}
    \caption{From left to right: ground-truth and $x_{IS}$ reconstruction of an image from the COULE data set.}
    \label{fig:gtEllissi}
\end{figure}

As visible from Figure \ref{fig:gtMAyo}, real full-dose medical images still present little noise and slightly visible streaking artifacts. This makes quite difficult to compute reliable metrics on the reconstructions. 
To fully exploit the full-reference image quality assessment metrics and validate our experiments, we create synthetic images and build few-view CT simulations, whose results are analysed in \ref{ssec:ResultsEllissi}.
In particular, the Constrasted Overlapping Uniform Lines and Ellipses (COULE) data set contains 430  sparse-gradient gray-scale images of $256 \times 256$ resolution with many overlying objects, varying in size and contrast with respect to the background.
The left image of Figure \ref{fig:gtEllissi} shows one image of the data set as an example.
The whole data set is downloadable from \url{www.kaggle.com/loiboresearchgroup/coule-dataset}.

The simulations are computed as before. We consider now different RIS images, i.e.  $x_{RIS} = x^{(3)}$, $x_{RIS} = x^{(5)}$ and $x_{RIS} = x^{(10)}$ for the $P_{360, 360}$ setting. 
We also explore different geometries $P_{360, 180}$ and $P_{360, 60}$, by reducing the number acquired projections.
On the right of Figure  \ref{fig:gtEllissi} we depict the SGP  solution $x_{IS}$ at convergence: we underline that it is visually not distinguishable from $x_{GT}$.


\subsection{The Network architecture and its training  \label{subsec:training}} %

As CNN architecture we use the state-of-the-art ResUNet architecture.
We exploit the pooling/unpooling strategies characterizing the popular U-net architecture by \cite{ronneberger2015u}, because we need to approximate the complex image-to-image transformation $f$.
Moreover, the addition of residual connections between the different levels of image resolution helps maintaining stable the training process and it has already been successfully applied to few-views CT image enhancement in \cite{han2018framing, ye2018deep}.
The activation function is $ReLU$ for each layer, except for the latter where $\tanh$ is used instead. \\
With Mayo Clinic data set we select $N = 3306$ images  (corresponding to 10 patients) for the training phase and 357 images (one patient) for testing. 
For the COULE data set, we use $N = 400$ images for training and 30 for testing. \\
Given the data set $\{ (x_{RIS, i}, x_{IS, i}) \}_{i=1, \dots, N}$, we denote by $F_{\theta^*}(x_{RIS, i})$ the action of the neural network on the input $x_{RIS, i}$ and  by $x_{ING, i} = F_{\theta^*}(x_{RIS, i})$ the network output. 
We estimate by network training, where the loss function is set as $\ell(x_{IS,i}, x_{ING,i}) = || x_{ING,i} - x_{IS,i} ||_2^2$, the parameters $\theta^*$ such that
\begin{equation}\label{eq:training}
    \theta^* = \arg\min_\theta \frac{1}{N} \sum_{i=1}^N \ell(x_{IS,i}, x_{ING,i}).
\end{equation}

The training is performed by running Adam for 50 and 100 epochs for the Mayo and COULE datasets, respectively. The batch size is fixed equal to 8 in all the experiments (this is the largest batch size usable  in our Nvidia RTX A4000 GPU).
The step size for the optimization algorithm decreases with polynomial decay, going from $10^{-3}$ to $10^{-5}$. To increase the stability over the first iterations, we clip the gradient to 5. 


\subsection{Implementation notes \label{subsec:implementation}} %
As already remarked in section \ref{subsec:RIS}, we use SGP as the iterative solver. 
In each iteration, we use the class OpTomo of Astra toolbox \cite{astra_1,astra_2} to compute the matrix-vector products $Ax$ and $A^Ty$, where $A$ is the matrix in equation \eqref{eq:LinearSyst}.
We remark that OpTomo provides the exact transpose of the projector operator $A$ for 2D tomography. 
In our experiments, we set the smoothing parameter for the TV regularizer to $\beta = 10^{-3}$, the regularization parameter as $\lambda=1 \cdot 10^{-5}$ for the  Mayo data set and $\lambda=4 \cdot 10^{-5}$ for the COULE data set. 
All the SGP inner parameters are taken from \cite{coap2018}. 


\subsection{Methods for comparisons \label{subsec:competitor}} %

In our experiments, we compare the RISING reconstructions with respect to the image solutions $x_{IS}$ achieved by the SGP solver at its numerical convergence, both to evaluate the visual improvements and  to analyse the learnability of the iterations.\\
In addition, we also consider the images computed by first applying the same RIS step as before and then restoring the $x_{RIS}$ image with a data-driven post-processing step. 
Here, we apply the same ResUNet architecture and  training setting used for RISING, but we consider  $\{ (x_{RIS, i}, x_{GT, i})\}_{i=1, \dots, N}$ as training set, imposing the ground truth images as targets (according to the standard LPP philosophy).
We denote such reconstructions as  $x_{LPP}$.

\section{Experimental results and discussion \label{sec:results}} 

In this section, we report and discuss the representative numerical experiments performed using the proposed  workflow for tomographic image reconstruction from few-view data. 

\subsection{Results on the real medical images  \label{ssec:ResultsMayo} }
\begin{figure*}[ht]
\centering
   \includegraphics[height=3.5cm]{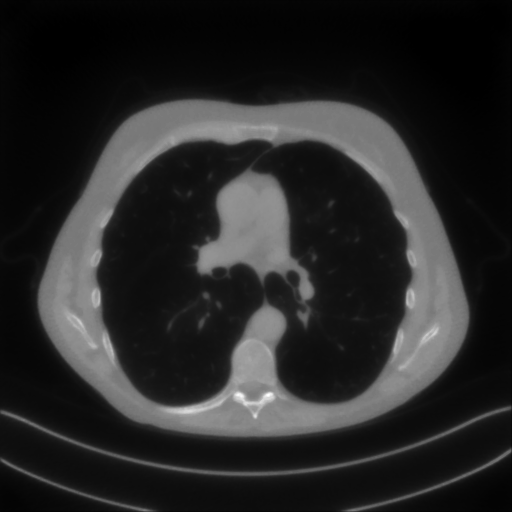}
 	   \includegraphics[trim= 190 122 130 160, clip, height=3.5cm]{imm/15A.png}
 	   \includegraphics[trim= 282 102 27 160, clip, height=3.5cm]{imm/15A.png} \qquad
    \includegraphics[height=3.5cm]{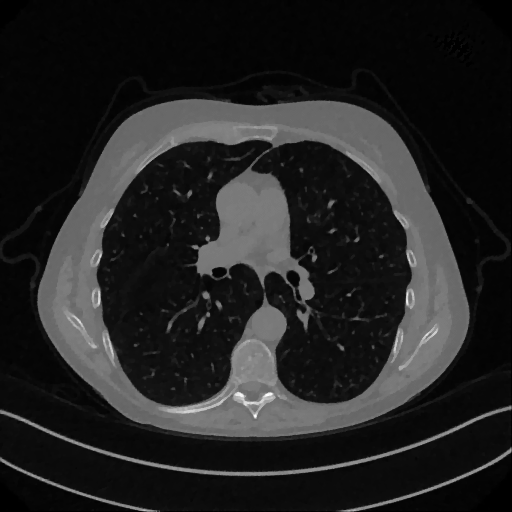}
       \includegraphics[trim= 190 122 130 160, clip, height=3.5cm]{imm/100A.png}
       \includegraphics[trim= 282 102 27 160, clip, height=3.5cm]{imm/100A.png}	\\ \vspace{5mm}
 	\includegraphics[height=3.5cm]{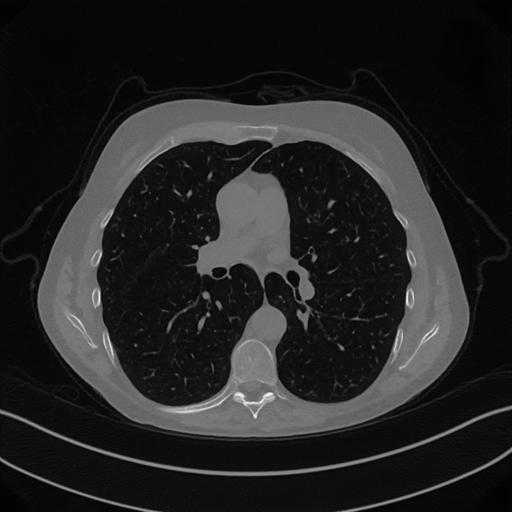}
 	   \includegraphics[trim= 190 122 130 160, clip, height=3.5cm]{imm/15AtoGT.png}
 	   \includegraphics[trim= 282 102 27 160, clip, height=3.5cm]{imm/15AtoGT.png} \qquad
    \includegraphics[height=3.5cm]{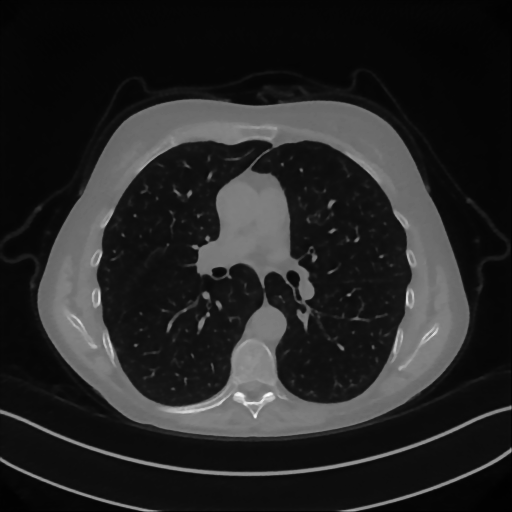}
 	    \includegraphics[trim= 190 122 130 160, clip, height=3.5cm]{imm/15Ato100A.png}
 	    \includegraphics[trim= 282 102 27 160, clip, height=3.5cm]{imm/15Ato100A.png}
\caption{Results on a test image from the Mayo data set, under the $P_{360, 360}$ CT protocol. Top-left: $x_{RIS}$; top-right: $x_{IS}$; bottom-left: $x_{LPP}$; bottom-right: $x_{ING}$. } 
\label{fig:testA}
\end{figure*}

\begin{figure*}[ht]
    \centering    \includegraphics[width=0.37\textwidth]{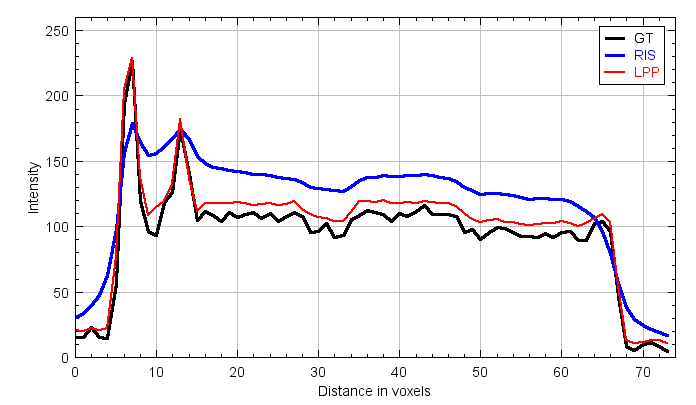}
    \includegraphics[width=0.37\textwidth]{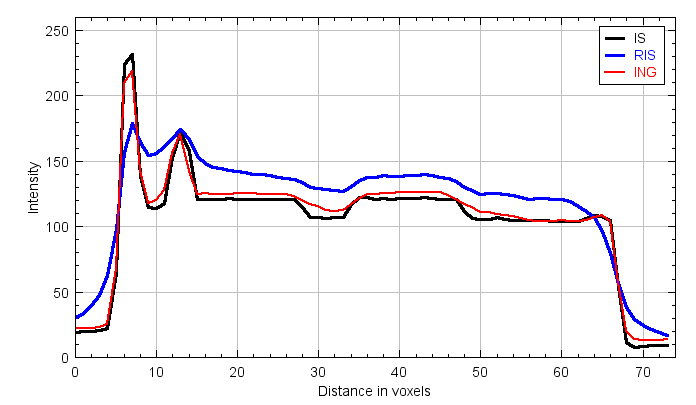}
    \caption{Intensity profiles taken on the horizontal red line depicted in \ref{fig:gtMAyo}, on the reconstructions in Figure \ref{fig:testA}. 
    Left: the black and red lines are the profile relative to $x_{GT}$ and $x_{LPP}$ respectively. Right: the black and red lines are the profile relative to $x_{IS}$ and $x_{ING}$ respectively. In both the plots, the blue line corresponds to the profile on the starting $x_{RIS}$ image.}
    \label{fig:Profili_A}
\end{figure*}

Here we present the results  of RISING applied to  the Mayo  data set introduced in paragraph \ref{sub:mayo}.
As previously mentioned, we consider two  sparse-view CT geometries, namely $P_{360, 360}$  and $P_{360, 180}$. 
In Figure \ref{fig:testA} we report the results for the  $P_{360, 360}$ protocol, achieved on one image of our test set.
The top-left image represents the $x_{RIS} = x^{(15)}$ reconstruction. Even if only a small number of iterations are performed, the main structures of the abdomen are visible; however, the image is still blurry.
In the $x_{IS}$ image,  shown in the upper-right corner, we notice that
 the TV regularizer has acted  to totally eliminate the artifacts and noise, improving the uniformity of the image in the inner structures. When compared to the ground-truth solution in Figure \ref{fig:gtMAyo}, the contours of the details in $x_{IS}$ appear slightly jagged, differently from $x_{GT}$ where they are neat but, usually, corrupted by artifacts. 
The bottom row of Figure \ref{fig:testA} shows the two $x_{LPP}$ and $x_{ING}$ images, respectively from left to right. 
It is evident that $x_{LPP}$ has retrieved many details but it presents noisy components, reflecting the features of its target image $x_{GT}$. 
Our solution $x_{ING}$ is less corrupted, since the low-contrast regions are correctly preserved and the noise is not visible. 
These observations are confirmed by Figure \ref{fig:Profili_A}, which plots the intensity profiles taken over the red line in the second crop (Figure \ref{fig:gtMAyo}).
In our approach (on the right) the CNN has accurately learnt the $f$ map of \eqref{eq:fComposition} and the $x_{ING}$ red profile mostly overlaps the black one.
On the contrary, the $x_{LPP}$ profile (on the left) is more distant from its target reference. 

\begin{figure*}[ht]
\centering
 	\includegraphics[height=3.5cm]{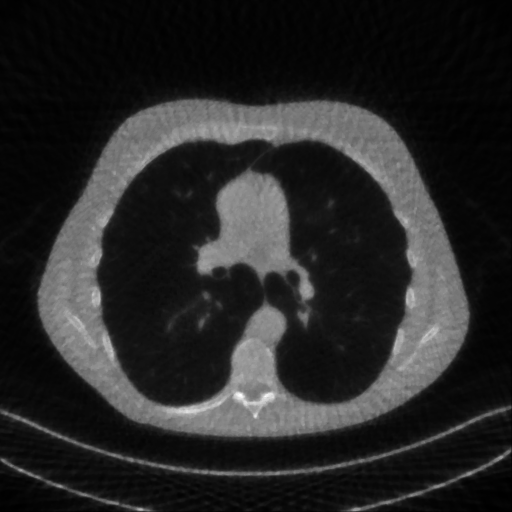}
 	   \includegraphics[trim= 190 122 130 160, clip, height=3.5cm]{imm/15B.png}
 	   \includegraphics[trim= 282 102 27 160, clip, height=3.5cm]{imm/15B.png} \qquad
 	\includegraphics[height=3.5cm]{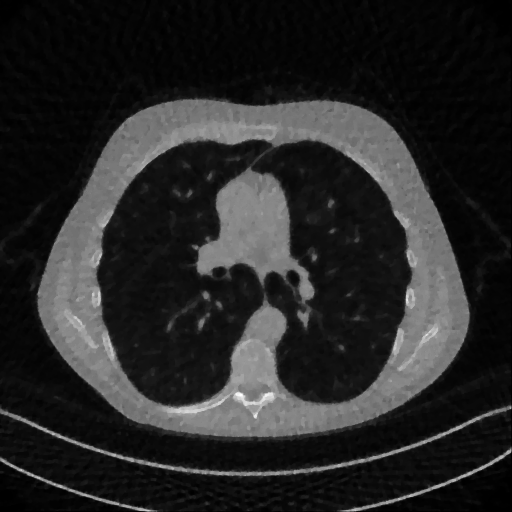}
        \includegraphics[trim= 190 122 130 160, clip, height=3.5cm]{imm/400B.png}
        \includegraphics[trim= 282 102 27 160, clip, height=3.5cm]{imm/400B.png}   \\ \vspace{5mm}
 	\includegraphics[height=3.5cm]{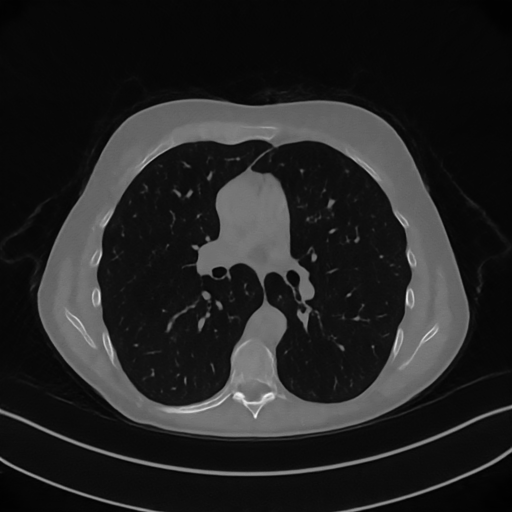}
 	    \includegraphics[trim= 190 122 130 160, clip, height=3.5cm]{imm/15BtoGT.png}
 	   \includegraphics[trim= 282 102 27 160, clip, height=3.5cm]{imm/15BtoGT.png} \qquad
    \includegraphics[height=3.5cm]{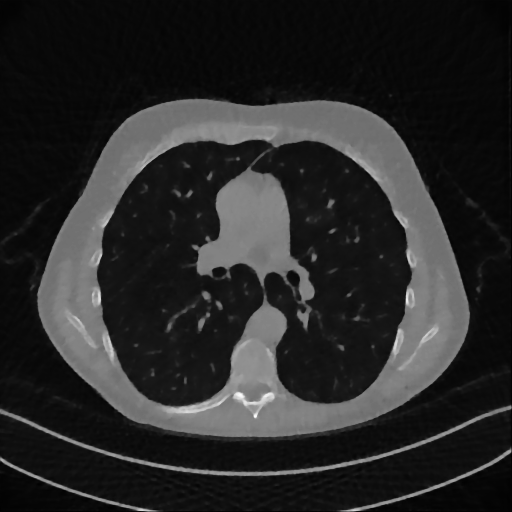}
 	    \includegraphics[trim= 190 122 130 160, clip, height=3.5cm]{imm/15Bto400B.png}
 	    \includegraphics[trim= 282 102 27 160, clip, height=3.5cm]{imm/15Bto400B.png}
\caption{Results on a test image from the Mayo data set, under the $P_{180, 60}$ CT protocol. Top-left: $x_{RIS}$; top-right: $x_{IS}$; bottom-left: $x_{LPP}$; bottom-right: $x_{ING}$. }
\label{fig:testB1}
\end{figure*}

\begin{figure*}[ht]
    \centering
    \includegraphics[width=0.37\textwidth]{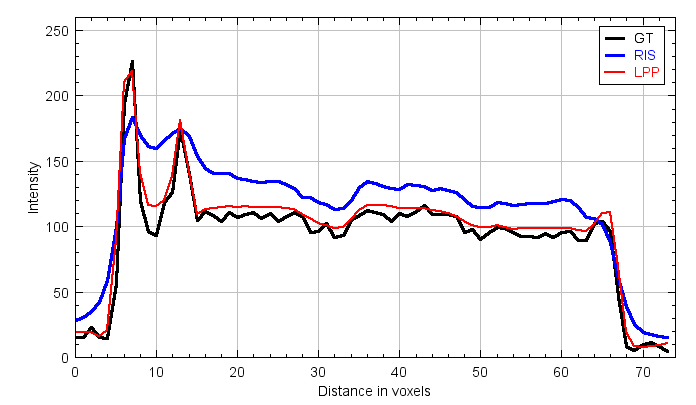}
    \includegraphics[width=0.37\textwidth]{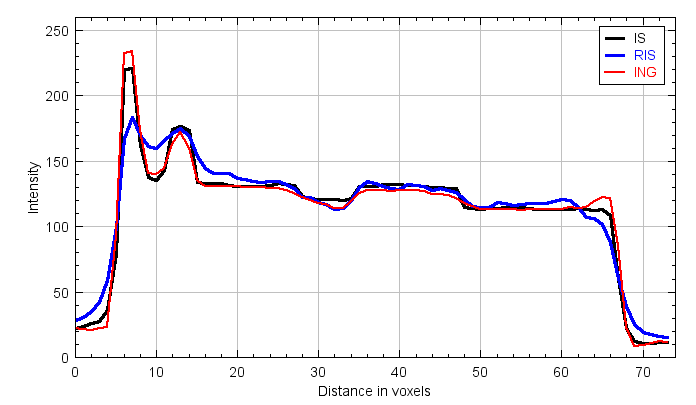}
    \caption{Intensity profiles taken on the horizontal red line depicted in \ref{fig:gtMAyo}, on the reconstructions in Figure \ref{fig:testB1}. 
    Left: the black and red lines are the profile relative to $x_{GT}$ and $x_{LPP}$ respectively. Right: the black and red lines are the profile relative to $x_{IS}$ and $x_{ING}$ respectively. In both the plots, the blue line corresponds to the profile on the starting $x_{RIS}$ image.}
    \label{fig:Profili_15B}
\end{figure*}

We now consider the $P_{180, 60}$ CT protocol whose results are reported in Figure \ref{fig:testB1}.
In this case, the tomographic  reconstruction is more challenging than in the previous experiment. The starting image $x_{RIS}$  has  evident streaking artifacts and blur and some details are lost, especially in the  first zoom.  The artifacts are reduced in the  $x_{IS}$ (top right image), where  some  details are recovered and the edges are quite neat.
The $x_{ING}$ image obtained with the proposed RISING (bottom right) is visually an excellent  reconstruction. It is  very similar to the $x_{LPP}$ image, whose training, we remark,  is  based on more informative target images. \\
In Figure \ref{fig:Profili_15B} we plot for these experiments, the same profiles of Figure \ref{fig:Profili_A}.
In the left graph, we analyse the performance of LPP approach: the $x_{LPP}$ reconstruction gets values quite close to the GT, but it does not fit well the target black line. 
In the right graph, we analyse the performance of RISING approach:  the $x_{ING}$ solution almost overlaps the target $x_{IS}$ profile, confirming that the network has correctly learned.

At last we underline that the solutions of the CS regularized model are very similar in case of $P_{360,360}$ and $P_{180,60}$ geometries.

\subsection{Results on synthetic images \label{ssec:ResultsEllissi} }

In this paragraph we report the results obtained with RISING on the synthetic COULE data set introduced in paragraph \ref{ssec:COULE}. Here, the  $x_{GT}$ images are digitally created, hence they are not corrupted by noise and artifacts and they can really be used as true images. 

\subsubsection{Metrics for image quality assessment\label{subsec:metrics}} 
To evaluate the efficiency and the accuracy of the reconstructed images, in the following we compute some full-reference image quality assessment metrics. 
For each considered image $x$, we consider the relative error (RE) of $x$ with respect to $x_{GT}$:
\begin{equation} \label{eq:RE}
    RE(x) = \frac{ \|x - x_{GT}\|_2^2 }{\|x_{GT}\|_2^2}
\end{equation} 
and the Root Mean Square Error (RMSE) between $x$ and the reference image $y$:
\begin{equation}
RMSE (x,y) = \sqrt{ \frac{\|x - y\|_2^2}{n^2} }
\label{eq:5}
\end{equation}
In the following, we set $y = x_{GT}$ unless otherwise stated.
We also compute the well-known Structural Similarity (SSIM) index   measuring the perceptual difference between two similar images $x$ and $y$, which is defined as in \cite{wang2004image}.
In our experiments, we compute the SSIM fixing $y =  x_{GT}$.

\subsubsection{Analysis for varying rapid reconstructions} 

\begin{figure*}[ht]
    \centering
    \includegraphics[width=0.2\textwidth]{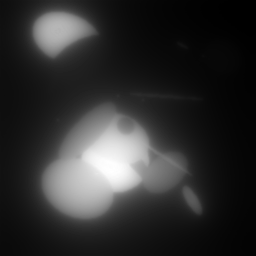}
    \includegraphics[width=0.2\textwidth]{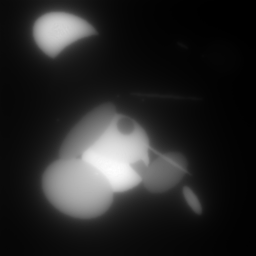}
    \includegraphics[width=0.2\textwidth]{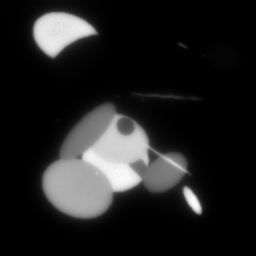} \\ \vspace{2mm}
    \includegraphics[width=0.2\textwidth]{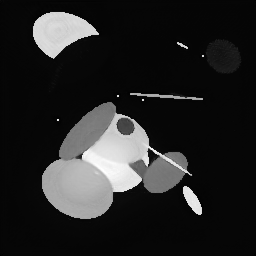}
    \includegraphics[width=0.2\textwidth]{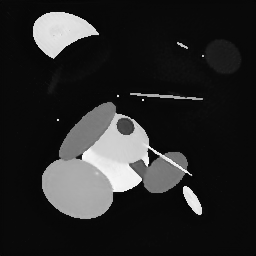}
    \includegraphics[width=0.2\textwidth]{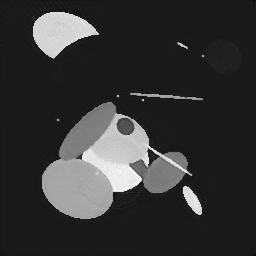}\\
    \caption{Results on a test image from the COULE synthetic data set, under the $P_{360, 360}$ CT protocol. First row, from left to right: $x_{RIS}$ with $K=3$, $K = 5$ and $K=10$ respectively; second row, from left to right: the corresponding $x_{ING}$.}
    \label{fig:Ellissi}
\end{figure*}

\begin{figure}[ht]
    \centering
    \includegraphics[trim= 0 0 0 40, clip, width=0.4\textwidth]{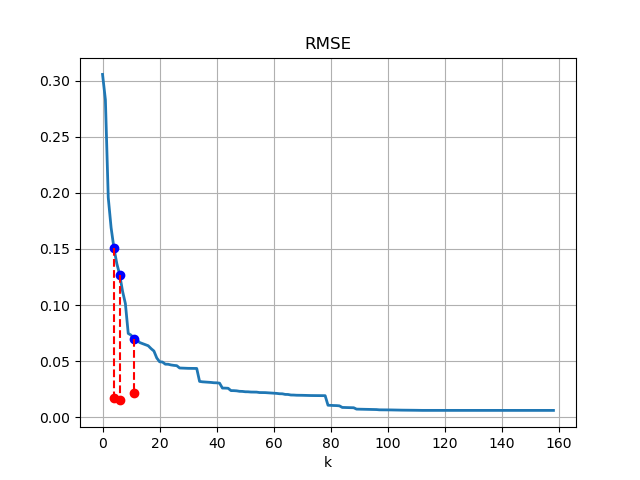} 
    \caption{Plot of the RMSE values over the performed iterations, for the reconstruction of an image from the COULE data set, in case of $P_{360, 360}$ geometry. The red dots are relative to the $x_{ING}$ images, achieved by the corresponding  $x^{(3)}$, $x^{(5)}$ and $x^{(10)}$ RIS images. }
    \label{fig:Ellissi_GraficoStarting}
\end{figure}

As first in-depth analysis, we set the geometry $P_{360, 360}$ and we compute the RIS early solutions for $K = \{3,5,10\}$. 
Figure \ref{fig:Ellissi} shows the $x_{RIS}$ starting images in the top row and the $x_{ING}$ final reconstructions in the bottom row (relative to the example image of Figure \ref{fig:gtEllissi}). They show that the network in the RIS step is able to almost perfectly learn from all three input images.
Figure \ref{fig:Ellissi_GraficoStarting} shows in blue the evolution of the RMSE through the iterations of the SGP for the considered image. In this case $K^*=58$, i.e. $x_{IS}=x^{(158)}$. The red dots represent the value of the errors of the $x_{ING}$ final reconstructions, for $K = \{3, 5, 10\}$. We observe that  the final $x_{ING}$ images are very close to the target $x_{IS}$, highlighting that the number $K$ of starting iterations seems not to notably influence the final results.  
In Table \ref{tab:Ellissi_Starting} we report the analysis on the whole test set, by means of the metrics introduced in the previous paragraph. We observe that the values relative to the  $x_{ING}$ solutions are all very similar, independently on $K$. Moreover,  the metrics for the $x_{IS}$ solution show that the image computed by the SGP method well approximates the true one and confirms that the proposed model-based reconstruction algorithm is very efficient.

\begin{table*}
\caption{Mean and standard deviation values of the quality metrics, evaluated on the COULE test set for different RIS reconstructions.}
    \label{tab:Ellissi_Starting}
  \begin{center}
 \begin{adjustbox}{max width=\textwidth}
\begin{tabular}{ll|cccc}
 & &  K=3 & K=5 & K=10  &    convergence \\
\thickhline
\multirow{3}{1cm}{RE}   & $x_{RIS}$ & 0.6152 $\pm$ 0.0366  & 0.4996 $\pm$ 0.0548 &  0.2700 $\pm$ 0.0239 & - \\
                        & $x_{ING}$ & 0.0681 $\pm$ 0.0098  & 0.0784 $\pm$ 0.0844  & 0.0781 $\pm$ 0.0122    & - \\
                        & $x_{IS}$ & -  & - &  - & 0.0207 $\pm$ 0.0051 \\
\thickhline
\multirow{3}{1cm}{RMSE} & $x_{RIS}$  & 0.1768 $\pm$ 0.0238  & 0.1432 $\pm$ 0.0218  & 0.0775 $\pm$ 0.0111  & -  \\
                        & $x_{ING}$  & 0.0196 $\pm$ 0.0042  & 0.0233 $\pm$ 0.0278  & 0.0225 $\pm$ 0.0053     & - \\
                        & $x_{IS}$  & -  & - & - & 0.0060 $\pm$ 0.0017   \\
\thickhline
\multirow{2}{1cm}{SSIM} & $x_{RIS}$  & 0.2264 $\pm$ 0.0262  & 0.3139 $\pm$ 0.0919  & 0.7844 $\pm$ 0.0481    & -  \\
                        & $x_{ING}$  & 0.9630 $\pm$ 0.0133  & 0.9267 $\pm$ 0.0389  & 0.9674 $\pm$ 0.0140     & - \\       
                        & $x_{IS}$  & - & -  & -   & 0.9980 $\pm$ 0.0010  \\
\thickhline
\end{tabular}
\end{adjustbox}
  \end{center}
\end{table*}

\subsubsection{Analysis for increasing sparsity of the CT subsampling} 
To understand the behavior of our approach at increasing sparsity in the CT protocol, we also test RISING at different geometric settings such as $P_{360,360}$, $P_{360,180}$ and $P_{360, 60}$. The number of iterations  used to generate $x_{RIS}$ is $K = 10$ and the $K^*$ iterations needed for convergence is for all the test images in $[150,300]$.
In Table \ref{tab:Ellissi_NTheta} we report the quality indexes evaluated on the test set in terms of mean and standard deviation.
We first interestingly observe that the values of all the $x_{RIS}$ outputs are very similar, independently from the  geometry. 
As before, the $x_{IS}$ solutions have excellent metrics, justifying the use of the RISING approach even in the hardest case with only 60 angles. This shows that the considered CS model \eqref{eq:7} properly describes the reconstruction process independently of the sparsity of the geometry (at least for the considered ones).
Concerning  final reconstruction $x_{ING}$, we see that halving the number of angles in $P_{360,180}$ does not affect the mean values of the metrics; when the geometry is very sparse in $P_{360,60}$ the errors slightly increase. However, the values of the standard deviations are small and very similar in all the tests, showing that the network has a stable behaviour.

\begin{table*}
  \caption{Mean and standard deviation values of the quality metrics, evaluated on the COULE test set for different geometries.}
    \label{tab:Ellissi_NTheta}
  \begin{center}
 \begin{adjustbox}{max width=\textwidth}
\begin{tabular}{ll|ccc}
                &  &  $P_{360,360}$ & $P_{360,180}$ & $P_{360, 60}$ \\
\thickhline
\multirow{3}{1cm}{RE}   & $x_{RIS}$  & 0.2700 $\pm$ 0.0239  & 0.2706 $\pm$ 0.0240  & 0.2834 $\pm$ 0.0282  \\
                        & $x_{ING}$  & 0.0781 $\pm$ 0.0122  & 0.0676 $\pm$ 0.0123  & 0.1113 $\pm$ 0.0223    \\
                        & $x_{IS}$   & 0.0207 $\pm$ 0.0051  & 0.0302 $\pm$ 0.0077  & 0.0668 $\pm$ 0.0147    \\
\thickhline
\multirow{3}{1cm}{RMSE} & $x_{RIS}$  & 0.0775 $\pm$ 0.0111  & 0.0777 $\pm$ 0.0111  & 0.0814 $\pm$ 0.0123    \\
                        & $x_{ING}$  & 0.0225 $\pm$ 0.0053  & 0.0198 $\pm$ 0.0062  & 0.0318 $\pm$ 0.0062     \\
                        & $x_{IS}$   & 0.0060 $\pm$ 0.0017  & 0.0088 $\pm$ 0.0026  & 0.0194 $\pm$ 0.0059    \\
\thickhline
\multirow{3}{1cm}{SSIM} & $x_{RIS}$  & 0.7844 $\pm$ 0.0481  & 0.7831 $\pm$ 0.0485  & 0.7472 $\pm$ 0.0643    \\
                        & $x_{ING}$  & 0.9674 $\pm$ 0.0140  & 0.9741 $\pm$ 0.0115  & 0.9493 $\pm$ 0.0117     \\  
                        & $x_{IS}$   & 0.9980 $\pm$ 0.0010  & 0.9951 $\pm$ 0.0030  & 0.9753 $\pm$ 0.0141    \\                  
\thickhline
\end{tabular}
\end{adjustbox}
  \end{center}
\end{table*}

\subsubsection{Empirical analysis of learnability} 
\begin{figure}[ht]
    \centering
    \includegraphics[width=0.4\textwidth]{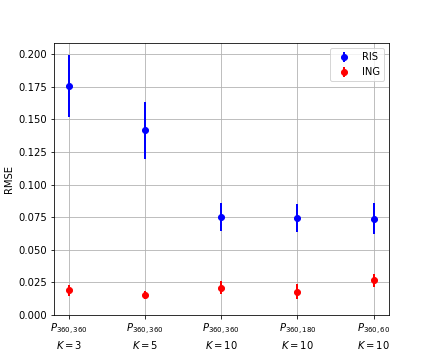}
    \caption{Error plots of the RMSE values, relative to the $x_{IS}$ target images, for the $x_{RIS}$ (in blue) and $x_{ING}$ (in red) reconstructions of the COULE images.}
    \label{fig:Learnability}
\end{figure}
In Figure \ref{fig:Learnability} we visualize the error plots of the RMSE values, evaluated as in \eqref{eq:5} with $y = x_{IS}$, relative to both the $x_{RIS}$ and $x_{ING}$ images, for all the previously considered settings.
We observe that the mean values for the $x_{ING}$ solutions are very small in all the cases. This empirically proves that the neural network has accurately learnt  the function $f_K$ in \eqref{eq:fComposition}, regardless $K$ and the geometry considered.
The three very similar RMSE mean values of the $x_{RIS}$ images (for $K=10$ and different geometries) indicate that, in these experiments, the error decreasing along the SGP iterations  does not depend on the sinogram sparsity.\\
Moreover, we notice that the standard deviation of the RIS images errors  are quite high, even if they decrease as $K$ increases. On the contrary, in all the considered experiments, the standard deviations  of the ING solutions are extremely small, denoting the stability of the proposed framework when applied to different images.

\section{Conclusions \label{sec:concl}}   

We have proposed a novel framework, called RISING, for the reconstruction  of a CT image from few-views.  RISING uses a Convolutional Neural Network lo learn the final iterations of an iterative   algorithms solving a gradient-exploiting regularization problem to reconstruct the image from sub-sampled data.\\
Hence, in RISING the network is trained on a data set obtained by computing the solution of the iterative solver from few-view projections, whereas
 in the usual LPP schemes the network is trained on {\em ground truth} images as targets. The numerical experiments performed 
 on a public data set of real abdomen images, the  reconstructions are visually very accurate, even in a very sparse geometry with 60 angles in $[0,180]$ degrees. Surprisingly, in this last case, the RISING network enhances the target image. 
On a simulated data set of gradient-sparse images, by means of error measures we numerically showed that the RISING framework computes an image very close to the ground truth, which is the desirable solution among the infinite possible of the inverse problems.
It is also noticeable that the errors do not depend neither on the iteration considered to get the coarse input image for the network nor on the acquisition geometry.

Finally, we underline that the RISING framework is very flexible, since it can be applied by considering any variational model for computing the inverse problem solution (not limiting to \eqref{eq:7}) and/or using other iterative solvers than SGP. In these cases, simply a different function $f_K$ is learned from the network. We also intend to study more theoretical properties, such as some conditions on the iterative methods and on the neural network which guarantee the convergence of the framework solution to the GT image.






\bibliographystyle{unsrturl} 
\bibliography{biblio_NNforCT}

\begin{thebibliography}{10}

\bibitem{mueller2012linear}
Jennifer~L Mueller and Samuli Siltanen.
\newblock {\em Linear and nonlinear inverse problems with practical
  applications}.
\newblock SIAM, 2012.

\bibitem{tuy1983inversion}
Heang~K Tuy.
\newblock An inversion formula for cone-beam reconstruction.
\newblock {\em SIAM Journal on Applied Mathematics}, 43(3):546--552, 1983.

\bibitem{donoho2006compressed}
David~L Donoho.
\newblock Compressed sensing.
\newblock {\em IEEE Transactions on information theory}, 52(4):1289--1306,
  2006.

\bibitem{GraffSidky2016}
C.~Graff and E.~Sidky.
\newblock Compressive sensing in medical imaging.
\newblock {\em Appl. Opt.}, 54(8):C23--C44, 2015.

\bibitem{wang2018image}
Ge~Wang, Jong~Chu Ye, Klaus Mueller, and Jeffrey~A Fessler.
\newblock Image reconstruction is a new frontier of machine learning.
\newblock {\em IEEE transactions on medical imaging}, 37(6):1289--1296, 2018.

\bibitem{han2016deep}
Yo~Seob Han, Jaejun Yoo, and Jong~Chul Ye.
\newblock Deep residual learning for compressed sensing ct reconstruction via
  persistent homology analysis.
\newblock {\em arXiv preprint arXiv:1611.06391}, 2016.

\bibitem{han2018framing}
Yoseob Han and Jong~Chul Ye.
\newblock Framing u-net via deep convolutional framelets: Application to
  sparse-view ct.
\newblock {\em IEEE transactions on medical imaging}, 37(6):1418--1429, 2018.

\bibitem{sidky2020cnns}
Emil~Y Sidky, Iris Lorente, Jovan~G Brankov, and Xiaochuan Pan.
\newblock Do cnns solve the ct inverse problem?
\newblock {\em IEEE Transactions on Biomedical Engineering}, 68(6):1799--1810,
  2020.

\bibitem{Scirep20}
R~Cavicchioli, J~Hu, E~Loli~Piccolomini, E~Morotti, and L~Zanni.
\newblock A first-order primal-dual algorithm for convex problems with
  applications to imaging. {GPU} acceleration of a model-based iterative method
  for digital breast tomosynthesis.
\newblock {\em Scientific Reports}, 10(1):120--145, 2020.
\newblock \href {https://doi.org/https://doi.org/10.1038/s41598-019-56920-y}
  {\path{doi:https://doi.org/10.1038/s41598-019-56920-y}}.

\bibitem{piccolomini2021model}
Elena~Loli Piccolomini and Elena Morotti.
\newblock A model-based optimization framework for iterative digital breast
  tomosynthesis image reconstruction.
\newblock {\em Journal of Imaging}, 7(2):36, 2021.

\bibitem{kak2001principles}
Avinash~C Kak and Malcolm Slaney.
\newblock {\em Principles of computerized tomographic imaging}.
\newblock SIAM, 2001.

\bibitem{candes2006robust}
Emmanuel~J Cand{\`e}s, Justin Romberg, and Terence Tao.
\newblock Robust uncertainty principles: Exact signal reconstruction from
  highly incomplete frequency information.
\newblock {\em IEEE Transactions on information theory}, 52(2):489--509, 2006.

\bibitem{piccolomini2016fast}
E~Loli Piccolomini and E~Morotti.
\newblock A fast total variation-based iterative algorithm for digital breast
  tomosynthesis image reconstruction.
\newblock {\em Journal of Algorithms \& Computational Technology},
  10(4):277--289, 2016.

\bibitem{xu2016accelerated}
Qiaofeng Xu, Deshan Yang, Jun Tan, Alex Sawatzky, and Mark~A Anastasio.
\newblock Accelerated fast iterative shrinkage thresholding algorithms for
  sparsity-regularized cone-beam ct image reconstruction.
\newblock {\em Medical physics}, 43(4):1849--1872, 2016.

\bibitem{sidky2014constrained}
Emil~Y Sidky, Rick Chartrand, John~M Boone, and Xiaochuan Pan.
\newblock Constrained {T}p{V} minimization for enhanced exploitation of
  gradient sparsity: Application to ct image reconstruction.
\newblock {\em IEEE journal of translational engineering in health and
  medicine}, 2:1--18, 2014.

\bibitem{purisha2017controlled}
Zenith Purisha, Juho Rimpel{\"a}inen, Tatiana Bubba, and Samuli Siltanen.
\newblock Controlled wavelet domain sparsity for x-ray tomography.
\newblock {\em Measurement Science and Technology}, 29(1):014002, 2017.

\bibitem{monga2021algorithm}
Vishal Monga, Yuelong Li, and Yonina~C Eldar.
\newblock Algorithm unrolling: Interpretable, efficient deep learning for
  signal and image processing.
\newblock {\em IEEE Signal Processing Magazine}, 38(2):18--44, 2021.

\bibitem{bengio_book}
Ian~J. Goodfellow, Yoshua Bengio, and Aaron Courville.
\newblock {\em Deep Learning}.
\newblock MIT Press, Cambridge, MA, USA, 2016.
\newblock \url{http://www.deeplearningbook.org}.

\bibitem{adler2017solving}
Jonas Adler and Ozan {\"O}ktem.
\newblock Solving ill-posed inverse problems using iterative deep neural
  networks.
\newblock {\em Inverse Problems}, 33(12):124007, 2017.

\bibitem{adler2018learned}
Jonas Adler and Ozan {\"O}ktem.
\newblock Learned primal-dual reconstruction.
\newblock {\em IEEE transactions on medical imaging}, 37(6):1322--1332, 2018.

\bibitem{gupta2018cnn}
Harshit Gupta, Kyong~Hwan Jin, Ha~Q Nguyen, Michael~T McCann, and Michael
  Unser.
\newblock Cnn-based projected gradient descent for consistent ct image
  reconstruction.
\newblock {\em IEEE transactions on medical imaging}, 37(6):1440--1453, 2018.

\bibitem{fista_net}
Jinxi Xiang, Yonggui Dong, and Yunjie Yang.
\newblock Fista-net: Learning a fast iterative shrinkage thresholding network
  for inverse problems in imaging.
\newblock {\em IEEE Transactions on Medical Imaging}, 40(5):1329--1339, 2021.
\newblock \href {https://doi.org/10.1109/TMI.2021.3054167}
  {\path{doi:10.1109/TMI.2021.3054167}}.

\bibitem{zhang2020metainv}
Haimiao Zhang, Baodong Liu, Hengyong Yu, and Bin Dong.
\newblock Metainv-net: Meta inversion network for sparse view ct image
  reconstruction.
\newblock {\em IEEE Transactions on Medical Imaging}, 40(2):621--634, 2020.

\bibitem{pelt2018improving}
Dani{\"e}l~M Pelt, Kees~Joost Batenburg, and James~A Sethian.
\newblock Improving tomographic reconstruction from limited data using
  mixed-scale dense convolutional neural networks.
\newblock {\em Journal of Imaging}, 4(11):128, 2018.

\bibitem{zhang2019dualres}
Tao Zhang, Hewei Gao, Yuxiang Xing, Zhiqiang Chen, and Li~Zhang.
\newblock Dualres-unet: Limited angle artifact reduction for computed
  tomography.
\newblock In {\em 2019 IEEE Nuclear Science Symposium and Medical Imaging
  Conference (NSS/MIC)}, pages 1--3. IEEE, 2019.

\bibitem{schnurr2019simulation}
Alena-Kathrin Schnurr, Khanlian Chung, Tom Russ, Lothar~R Schad, and Frank~G
  Z{\"o}llner.
\newblock Simulation-based deep artifact correction with convolutional neural
  networks for limited angle artifacts.
\newblock {\em Zeitschrift f{\"u}r Medizinische Physik}, 29(2):150--161, 2019.

\bibitem{urase2020simulation}
Yasuyo Urase, Mizuho Nishio, Yoshiko Ueno, Atsushi~K Kono, Keitaro Sofue,
  Tomonori Kanda, Takaki Maeda, Munenobu Nogami, Masatoshi Hori, and Takamichi
  Murakami.
\newblock Simulation study of low-dose sparse-sampling ct with deep
  learning-based reconstruction: usefulness for evaluation of ovarian cancer
  metastasis.
\newblock {\em Applied Sciences}, 10(13):4446, 2020.

\bibitem{morotti2021green}
Elena Morotti, Davide Evangelista, and Elena Loli~Piccolomini.
\newblock A green prospective for learned post-processing in sparse-view
  tomographic reconstruction.
\newblock {\em Journal of Imaging}, 7(8):139, 2021.

\bibitem{bubba2019learning}
Tatiana~A Bubba, Gitta Kutyniok, Matti Lassas, Maximilian Maerz, Wojciech
  Samek, Samuli Siltanen, and Vignesh Srinivasan.
\newblock Learning the invisible: a hybrid deep learning-shearlet framework for
  limited angle computed tomography.
\newblock {\em Inverse Problems}, 35(6):064002, 2019.

\bibitem{jiang2019augmentation}
Zhuoran Jiang, Yingxuan Chen, Yawei Zhang, Yun Ge, Fang-Fang Yin, and Lei Ren.
\newblock Augmentation of cbct reconstructed from under-sampled projections
  using deep learning.
\newblock {\em IEEE transactions on medical imaging}, 38(11):2705--2715, 2019.

\bibitem{RUDIN1992}
Leonid~I. Rudin, Stanley Osher, and Emad Fatemi.
\newblock Nonlinear total variation based noise removal algorithms.
\newblock {\em Physica D: Nonlinear Phenomena}, 60(1):259 -- 268, 1992.
\newblock URL:
  \url{http://www.sciencedirect.com/science/article/pii/016727899290242F},
  \href {https://doi.org/https://doi.org/10.1016/0167-2789(92)90242-F}
  {\path{doi:https://doi.org/10.1016/0167-2789(92)90242-F}}.

\bibitem{bonettini2008scaled}
Silvia Bonettini, Riccardo Zanella, and Luca Zanni.
\newblock A scaled gradient projection method for constrained image deblurring.
\newblock {\em Inverse problems}, 25(1):015002, 2008.

\bibitem{coap2018}
Elena Loli~Piccolomini, V.L. Coli, E.~Morotti, and L.~Zanni.
\newblock Reconstruction of 3{D X}-ray {CT} images from reduced sampling by a
  scaled gradient projection algorithm.
\newblock {\em Comp. Opt. Appl.}, 71:171--191, 2018.
\newblock \href {https://doi.org/https://doi.org/10.1007/s10589-017-9961-2}
  {\path{doi:https://doi.org/10.1007/s10589-017-9961-2}}.

\bibitem{BPR16}
S~Bonettini, F~Porta, and V~Ruggiero.
\newblock A variable metric inertial method for convex optimization.
\newblock {\em SIAM J. Sci. Comput}, 31(4):A2558--A2584, 2016.

\bibitem{mccollough2016tu}
C~McCollough.
\newblock Tu-fg-207a-04: Overview of the low dose ct grand challenge.
\newblock {\em Medical physics}, 43(6Part35):3759--3760, 2016.

\bibitem{ronneberger2015u}
Olaf Ronneberger, Philipp Fischer, and Thomas Brox.
\newblock U-net: Convolutional networks for biomedical image segmentation.
\newblock In {\em International Conference on Medical image computing and
  computer-assisted intervention}, pages 234--241. Springer, 2015.

\bibitem{ye2018deep}
Jong~Chul Ye, Yoseob Han, and Eunju Cha.
\newblock Deep convolutional framelets: A general deep learning framework for
  inverse problems.
\newblock {\em SIAM Journal on Imaging Sciences}, 11(2):991--1048, 2018.

\bibitem{astra_1}
Wim {van Aarle}, Willem~Jan Palenstijn, Jan {De Beenhouwer}, Thomas Altantzis,
  Sara Bals, K.~Joost Batenburg, and Jan Sijbers.
\newblock The astra toolbox: A platform for advanced algorithm development in
  electron tomography.
\newblock {\em Ultramicroscopy}, 157:35--47, 2015.
\newblock URL:
  \url{https://www.sciencedirect.com/science/article/pii/S0304399115001060},
  \href {https://doi.org/https://doi.org/10.1016/j.ultramic.2015.05.002}
  {\path{doi:https://doi.org/10.1016/j.ultramic.2015.05.002}}.

\bibitem{astra_2}
Wim van Aarle, Willem~Jan Palenstijn, Jeroen Cant, Eline Janssens, Folkert
  Bleichrodt, Andrei Dabravolski, Jan~De Beenhouwer, K.~Joost Batenburg, and
  Jan Sijbers.
\newblock Fast and flexible x-ray tomography using the astra toolbox.
\newblock {\em Opt. Express}, 24(22):25129--25147, Oct 2016.
\newblock URL:
  \url{http://www.osapublishing.org/oe/abstract.cfm?URI=oe-24-22-25129}, \href
  {https://doi.org/10.1364/OE.24.025129} {\path{doi:10.1364/OE.24.025129}}.

\bibitem{wang2004image}
Zhou Wang, Alan~C Bovik, Hamid~R Sheikh, and Eero~P Simoncelli.
\newblock Image quality assessment: from error visibility to structural
  similarity.
\newblock {\em IEEE transactions on image processing}, 13(4):600--612, 2004.

\end{thebibliography}

\end{document}